\documentclass{amsart}
\usepackage{amsfonts}
\usepackage{amssymb}
\begin{document}

{\bf
Convolution-multiplication identities for Tutte polynomials of matroids}

\vskip .33in
\ \ Joseph P. S. Kung

\ \ Department of Mathematics,

\ \ University of North Texas, Denton, TX 76203, U.S.A.

\vskip .10in
e-mail: kung@unt.edu
\vskip .33in


\newcommand\join{\vee}
\newcommand\meet{\wedge}
\newcommand\rk{\mathrm{rk}}
\newcommand\SC{\mathrm{SC}}
\newcommand\NC{\mathrm{NC}}
\newcommand\RR{\mathrm{R}}

{\it Abstract.}
We give a general multiplication-convolution identity for the multivariate and bivariate
rank generating polynomial of a matroid. The bivariate rank generating polynomial is transformable to and from
the Tutte polynomial by simple algebraic operations.  Several identities,
almost all already known in some form, are specialization of this identity.  Combinatorial 
or probabilistic interpretations are given for the specialized identities.

{\it Subject classification number.}  05B35  05C15  05C80


\vskip .33in \noindent
\begin{center}
1.  A convolution-multiplication identity
\end{center}

This paper originated in an attempt to understand the relationship between
Identities 7 and 8 in this paper.   One way to do this is to find a
hidden identity which contains both identities as special cases.  
Identity 1 and its specialization,
Identity 2, provide candidates for such an identity.
As the reader will see, the proofs of Identities 1 and 2, and perhaps the
identities themselves, are trivial.  Identities have a way of receding to
a {\it Zen} state of triviality.
In the other direction, we add meaning by giving combinatorial or probabilistic interpretations for many identities derived from the hidden identity.  

Let $M$ be a rank-$r$ matroid on the set $E$ with rank function $\rk$ and
 $L(M)$ be its lattice of closed sets or flats.
If $M$ has no loops (that is, elements of rank zero), its
{\sl characteristic polynomial} $\chi (M;\lambda)$
is the polynomial in the variable $\lambda$ defined by
$$
\chi (M;\lambda)
 = \sum_{X:\, X \in L(M)} \mu(\hat{0},X) \lambda^{r - \rk(X)},
$$
where $\mu$ is the M\"obius function in the lattice $L(M).$  If $M$ has a loop,
then $\chi (M;\lambda)$ is defined to be the zero polynomial.
The {\sl size-corank polynomial} $\SC (M;x,\lambda)$ is
the polynomial in two variables, $x$ and $\lambda$ defined by
$$
\SC(M;x,\lambda) = \sum_{A:\,A \subseteq E}  x^{|A|} \lambda^{r - \rk (A)}.
$$
The size-corank polynomial has a natural ``polarized'' multivariate generalization.
Let $\underline{x}$ be a labeled multiset
$\{x_e: e \in E\}$ of variables or numbers, one
for each element $e$ in the ground set $E$ of the matroid $M.$   If
$A \subseteq E,$ define $\underline{x}^A$ to be the product
$\prod_{e:\, e \in A} x_e.$  This notation is used analogously, so that,
for example, $\underline{-x} = \{-x_e:\, e \in E\}$ and
$\underline{(-x)}^A = \prod_{e:\, e \in A} (-x_e).$  In addition, if
$\{y_e:\, e \in E\}$ is another multiset, then
$\underline{xy}$ is the set $\{x_e y_e: e \in E\}$
and $\underline{(xy)}^A = \prod_{e:\, e \in A} x_e y_e.$
The {\sl subset-corank polynomial} $\mathbf{SC}(M; \underline{x},\lambda)$
is defined by
$$
\mathbf{SC}(M; \underline{x},\lambda) =
 \sum_{A:\,A \subseteq E}  \underline{x}^A \lambda^{r - \rk (A)}.
$$
The subset-corank polynomial specializes to the size-corank  polynomial
when we set all the variables $x_e$ to the same variable $x.$
Variants of the subset-corank polynomial, usually defined for graphs, have
been rediscovered many times.  The original discovery
is by Fortuin and Kasteleyn in statistical mechanics.  An {\it almost} complete
list can be made by merging the lists in the surveys of Farr \cite{Farr}
and Sokal \cite{Sokal}.
For matroids, the subset-corank polynomial is due to R.T. Tugger and it is sometimes named after her (see \cite{Twelve}).

The {\sl nullity-corank polynomial,} usually known as the {\sl rank
generating polynomial,} $\RR (M;x,\lambda)$ is defined by
$$
\RR (M;x,\lambda) = \sum_{A:\,A \subseteq E}  x^{|A| - \rk(A)}
\lambda^{r - \rk (A)}.
$$
The size-corank and nullity-corank polynomials are closely related: indeed,
$$
x^{-r} \SC(M;x, x\lambda) = \RR (M;x,\lambda)
$$
and
$$
\SC(M;x, \lambda) = x^r \RR (M;x,\lambda/x).
$$
Both polynomials specialize to the characteristic polynomial.  Specifically,
$$
\SC(M;-1,\lambda) = \chi (M;\lambda), \,\,\,\, 
\RR(M;-1,-\lambda)= (-1)^r \chi (M;\lambda).
\eqno(1)$$
The {\sl Tutte polynomial} is defined to be the polynomial
$\RR(M;x-1, \lambda-1).$  The three bivariate
 polynomials transform into each other by simple algebra.  (Note, however, that
 we may need to divide; thus, problems may arise when we wish
 to set a variable to $0,$)
  We shall choose
 the polynomial that gives the simplest or most general identities.  Thus, despite the title,
 we shall work mostly with subset-corank or size-corank polynomials.

We begin with a general multiplication-convolution identity.

\vskip 0.2in\noindent
{\bf Identity 1.}   Let $M$ be a matroid on the set $E,$  $\underline{x}$
and $\underline{y}$ be multisets of variables labeled by $E,$
and $\lambda$ and $\xi$ be variables.  Then
$$
\mathbf{SC} (M; \underline{xy}, \lambda \xi)
=
\sum_{T:\, T \subseteq E} \lambda^{r - \rk (T)} \underline{(-y)}^{T}
\mathbf{SC} (M|T;\underline{-x},\lambda)
\mathbf{SC} (M/T;\underline{y},\xi).
$$

{\it Proof.}  We use the fact that if $B$ and $A$ are subsets of $E$ such that
$B \subseteq A,$ then the sum
$$
\sum_{T:\, B \subseteq T \subseteq A} (-1)^{|T|-|B|}
$$
equals $0$ except in the case $A = B,$ when it equals $1.$   Then

\begin{eqnarray*}
&&
\mathbf{SC} (M;\underline{xy},\lambda \xi)
\\
&=&
\sum_{B,A:  \, B,A \subseteq E}
\underline{x}^{B} \lambda^{r - \rk (B)} \underline{y}^{A} \xi^{r - \rk (A)}
\left[
\sum_{T:\, B \subseteq T \subseteq A} (-1)^{|T|-|B|}
\right]
\\
&=&
\sum_{T:  \, T \subseteq E}
\lambda^{r - \rk (T)} \underline{(-y)}^{T}
\left[\sum_{B:\,B \subseteq T}  \underline{(-x)}^{B} \lambda^{r - \rk(B)}
\right]
\left[\sum_{A:\,T \subseteq A \subseteq E}  \underline{y}^{A \backslash T}
\xi^{(r - \rk(T)) - (\rk(A) - \rk(T))}
\right]
\\
&=&
\sum_{T:\, T \subseteq E} \lambda^{r - \rk (T)} \underline{(-y)}^{T}
\mathbf{SC} (M|T;\underline{-x},\lambda) \mathbf{SC} (M/T;\underline{y},\xi).
\end{eqnarray*}
\qed

\vskip 0.2in\noindent
Setting $x_e = x,$ we obtain the following specialization of Identity 1.

\vskip 0.2in\noindent
{\bf Identity 2.}  Let $x,y, \lambda,\xi$ be variables.  Then
$$
\SC (M;xy,\lambda \xi)
=
\sum_{T:\, T \subseteq E} \lambda^{r - \rk (T)} (-y)^{|T|}
\SC (M|T;-x,\lambda)
\SC (M/T;y,\xi).
$$

Next we state, without proof, the analog of Identity 2 for the
nullity-corank polynomial.

\vskip 0.2in\noindent
{\bf Identity 3.}
$$
\RR (M;xy,\lambda \xi)
=
\sum_{T:\, T \subseteq E} \lambda^{r - \rk (T)} (-y)^{|T| - \rk(T)}
\RR (M|T;-x,-\lambda)
\RR (M/T;y,\xi).
$$

\vskip 0.4in\noindent
The algebraic operation dual to convolution is comultiplication.  Thus, we 
can express the identities in this section in the language of coalgebras.  However, we will wait until we have more than just a formal theory.

\vskip 0.4in\noindent
***Proof of Identity 3.\footnote{
This proof is included so that one can easily checked the identity.  As indicated
in the text, it
will be removed in the final version.  

\begin{eqnarray*}
&&
\RR (M;xy,\lambda \xi)
\\
&=&
\sum_{B,A:  \, B,A \subseteq E}
x^{|B| - \rk (B)} \lambda^{r - \rk (B)} y^{|A|- \rk(A)} \xi^{r - \rk (A)}
\left[
\sum_{T:\, B \subseteq T \subseteq A} (-1)^{|T|-|B|}
\right]
\\
&=&
\sum_{T:  \, T \subseteq E}
\lambda^{r - \rk (T)} (-y)^{|T| - \rk(T)}
\left[\sum_{B:\,B \subseteq T}  (-x)^{|B|-\rk(B)} (-\lambda)^{\rk(T) - \rk(B)}
\right]
\\
&& \qquad
\cdot \left[\sum_{A:\,T \subseteq A \subseteq E}  y^{(|A| - |T|)- (\rk(A) - \rk(T))}
\xi^{(r - \rk(T)) - (\rk(A) - \rk(T))}
\right]
\\
&=&
\sum_{T:\, T \subseteq E}
\lambda^{r - \rk (T)} (-y)^{|T| - \rk(T)}
\RR (M|T;-x,-\lambda)  \RR (M/T;y,\xi).
\end{eqnarray*}
}

\vskip .33in \noindent
\begin{center}
2.  Weighted sums of polynomials
\end{center}

\vskip 0.2in\noindent
Identities 1, 2. and 3 specialize to several known identities.   We begin
with identities expressing the size-corank or subset-corank polynomial
as a weighted sum of other polynomials.

\vskip 0.2in\noindent
{\bf Identity 4.}
\begin{eqnarray*}
\mathbf{SC} (M; \underline{xy},\xi) =
\sum_{T:\, T \subseteq E} \underline{(x+1)}^T \underline{y}^T
\SC (M/T;\underline{-y},\xi).
\end{eqnarray*}

{\it Proof.}  Set $y_e = -y_e,$ $\lambda = 1,$ $x_e = -x_e,$
leaving $\xi$ unchanged, in Identity 1 to obtain
\begin{eqnarray*}
\mathbf{SC} (M;\underline{xy},\xi) &=&
\sum_{T:\, T \subseteq E} \underline{y}^T \mathbf{SC} (M|T; \underline{x},1)
\mathbf{SC} (M/T;\underline{-y},\xi)
\\
&=&
\sum_{T:\, T \subseteq E} \underline{(x+1)}^{T} \underline{y}^T\mathbf{SC} (M/T;\underline{-y},\xi).
\end{eqnarray*}
In the last step in the derivation, we used a multivariate version
of the binomial identity:
\begin{eqnarray*}
\mathbf{SC} (M|T; \underline{x},1) &=&
\sum_{B:\, B \subseteq T} \,\, \prod_{b:\, b \in B} x_b
\\
&=&
\prod_{e:\, e \in T} (x_e + 1)
= \underline{(x+1)}^T.
\end{eqnarray*}
\qed

\vskip 0.2in\noindent
Identity 4 specializes to more familiar identities.

\vskip 0.2in\noindent
{\bf Identity 5.}
\begin{eqnarray*}
\mathbf{SC} (M; \underline{x},\xi) &=&
\sum_{X:\, X \in L(M)} \underline{(x+1)}^X \chi (M/X;\xi),
\\
\SC (M;x,\xi) &=& \sum_{X:\, X \in L(M)} (x+1)^{|X|} \chi (M/X;\xi),
\end{eqnarray*}
where the sums range over all closed sets $X$ in the lattice $L(M).$

{\it Proof.}
If $T$ is not a closed set, then
$\SC (M/T;-1,\xi) = 0$
and if $T$ is a closed set,  then $ \SC (M/T; -1,\xi) = \chi (M/T;\xi). $  Thus the
range of the sums can be restricted to closed sets.
\qed

\vskip 0.2in\noindent
The bivariate form of Identity 5 is a fundamental identity of Tutte \cite{Tutte}.  Tutte found it for graphs and
Crapo \cite{Crapo} extended it to matroids. 
The multivariate form appeared in \cite{Twelve}.  It has the following
interpretation.  
The subset-corank polynomial encodes the rank function of the matroid $M$ in
the sense that the rank of a set $A$ is $r-d,$ where $d$ is the degree of
$\lambda$ in the monomial $\underline{x}^A \lambda^d$ in
$\mathbf{SC} (M;\underline{x},\lambda).$  On the other hand,
$\mathbf{SC} (M; \underline{x-1},\xi)$ encodes
the collection of closed sets of $M$
in the following way:  a set $T$ is closed
in $M$ if and only if the monomial $c\underline{x}^T$ occurs
in $\mathbf{SC} (M; \underline{x-1},\xi)$ with a nonzero coefficient $c.$  (The nonzero
coefficient is $\chi(M/T;\lambda).$)   Thus, Identity 5 gives
an ``algebraic transformation'' of the rank description of a matroid to its
closed set description.

\vskip .33in \noindent
\begin{center}
3.  Random matroids
\end{center}

\vskip .20in \noindent
We next consider identities which give interpretations of size-corank or subset-corank
polynomials as expected values of enumerative invariants of
a random submatroid.

We begin with an identity which is an ``order dual'' of Identity 5.
Setting $x_e=1,$ $\xi=1,$ $y_e = -y_e,$ and leaving $\lambda$ unchanged
in Identity 1, we obtain
\begin{eqnarray*}
\mathbf{SC} (M; \underline{-y},\lambda)
=
\sum_{T:\, T \subseteq E} \lambda^{r - \rk (T)}
y^{T}
\mathbf{SC} (M|T;\underline{-1},\lambda) \mathbf{SC} (M/T;\underline{-y},1).
\end{eqnarray*}
Using the fact that
$\mathbf{SC}(M/T;\underline{-y},1) = \underline{1-y}^{E \backslash T},$
 we obtain the following identity.

\vskip 0.2in\noindent
{\bf Identity 6.}
%
%
$$
\mathbf{SC} (M; \underline{-y}, \lambda)
=
\sum_{T:\, T \subseteq E} \lambda^{r - \rk (T)}
\chi (M|T;\lambda)
\underline{y}^T \underline{(1 - y)}^{E \backslash T}.
$$

\vskip 0.2in \noindent
The next identity, obtained by setting $x_e = - x_e,$ $y_e = -y_e,$ $\xi = 1,$
leaving $\lambda$ unchanged, is a generalization of Identity 6.

\vskip 0.2in\noindent
{\bf Identity 7.}
%
$$
\mathbf{SC} (M; \underline{xy}, \lambda)
=
\sum_{T:\, T \subseteq E} \lambda^{r - \rk (T)}
\mathbf{SC} (M|T; \underline{x},\lambda)
\underline{y}^T \underline{(1 - y)}^{E \backslash T}.
$$

\vskip 0.2in \noindent
Bivariate versions of Identities 6 and 7,
stated for graphs and given interpretations in terms of random graphs,
are known.  See Welsh \cite{Welsh}.  Specifically,
Identity 6 is related to
an identity in Vertigan's
Oxford thesis \cite{Vert} and Welsh has given a proof
using random subgraphs in \cite{Welsh}.   A version of Identity 7 was
found by Grimmett \cite{Grim} (see also \cite{Welsh}).

The interpretations by random subgraphs (given in \cite{Welsh})
generalize easily to interpretations by random submatroids.
Let $M$ be a matroid on the set $E.$  We generate a {\sl random subset}
and hence, a {\sl random submatroid}, of
$M$ by deleting each element $e$ in $E$ independently and at random
with probability $1 - p_e.$
We need two somewhat artificial definitions, designed to make the
interpretations work.
Let $N$ be a rank-$s$ submatroid of the rank-$r$ matroid $M.$  Then
the {\sl normalized} characteristic
polynomial $\chi^{\dagger}(N;\lambda)$ to be $\lambda^{r-s}\chi(N;\lambda).$
Similarly,
the {\sl normalized} size-corank polynomial $\mathrm{SC}^{\dagger}(N;x,\lambda)$
is defined to be $\lambda^{r - s} \mathrm{SC} (N;x,\lambda).$
Since the probability that the subset $T$ is
chosen is $\underline{p}^T \underline{(1-p)}^{E \backslash T},$
it is clear that the expected value
of the normalized characteristic (respectively, normalized
size-corank polynomial) of a random submatroid of $M$
equal $\mathbf{SC} (M;\underline{-p},\lambda)$ (respectively,
$\mathbf{SC} (M;\underline{px},\lambda)$).

In analogy with random graphs, one can develop a
theory of random sets of vectors or points in a finite {\sl ambient} space.
This was done in Kelly and Oxley \cite{KOxley}.
There are two choices for the ambient space:
the (affine) vector space $\mathrm{GF}(q)^d$ or the
projective space $\mathrm{PG}(d-1,q)$ over $\mathrm{GF}(q).$  Here,
$q$ is a prime power and $\mathrm{GF}(q)$ is the finite field of order $q.$
As a matroid, the projective space $\mathrm{PG}(d-1,q)$ is a simplification of
$\mathrm{GF}(q)^d.$ In particular, the two ambient spaces
have the same lattice $L(d,q)$ of flats.
The two ambient spaces give the same theory of random sets, more or less.
Since probability theorists usually work with random matrices and vectors,
we shall work with $\mathrm{GF}(q)^d.$  We define a {\sl random set $V(d,q,p)$
of vectors} to be a random submatroid of $\mathrm{GF}(q)^d$ in which each
element is chosen with the same probability $p.$

\vskip 0.2in \noindent
{\bf Lemma.}
\begin{eqnarray*}
\SC (\mathrm{GF}(q)^d; x, \lambda) =
\sum_{j=0}^d   {\binom {d}{j}}_q
(\lambda - 1) (\lambda - q) \cdots (\lambda - q^{d-i}) (x+1)^{q^j}
\end{eqnarray*}
%
%
$$
\sum_{d=0}^{\infty}  \SC (\mathrm{GF}(q)^d;x,\lambda) \frac {z^d}{[d!]_q}
=
\left[ \prod_{d=0}^{\infty} \, \frac { 1 + zq^d} { 1 + \lambda z q^d} \right]
\left[
\sum_{d=0}^{\infty}  (x+1)^{q^d} \frac {z^d}{[d!]_q}
\right],
$$
where 
$$
[d!]_q = (1-q)(1-q^2) \cdots (1-q^d) \,\,\mathrm{and} \,\,
{\binom {d}{j}}_q = \frac {[d!]_q} {[j!]_q [(d-j)!]_q}.
$$

{\it Proof.}  The formula for $\SC (\mathrm{GF}(q)^d;x,\lambda)$ follows
easily from known counting formulas in finite vector spaces and Identity 5.
It is also a special case of formulas in Mphako \cite{Mp}.  The generating
function of $\SC (\mathrm{GF}(q)^d;x,\lambda)$ factors into the product
$$
\left[
\sum_{d=0}^{\infty}
(\lambda - 1)  (\lambda - q) \cdots (\lambda - q^d) \frac {z^d}{[d!]_q}
\right]
\left[
\sum_{d=0}^{\infty}  (x + 1)^{q^d} \frac {z^d}{[d!]_q}
\right].
$$
The formula now follows from the $q$-binomial theorem, in the version stated
in \cite{GasRah},  Section 1.3.
\qed

\vskip 0.2in \noindent
Similar formulas exist for $\SC(\mathrm{PG}(d-1,q);x,\lambda).$  Simply replace the exponent $q^j$
in $x^{q^j}$ by $ q^{j-1}+q^{j-2}+\cdots+q+1.$

We next give two typical results about
expected values of matroid invariants of random vectors.
We shall need to assume some knowledge of critical problems
(see, for example, \cite{BryOx, KungC}).

\vskip 0.2in \noindent
{\bf Theorem.}
(a)
The expected number $\mathrm{D}(d,p,q,s)$ of $s$-tuples of linear functionals
on $\mathrm{GF}(q)^d$
distinguishing a random set $V(d,p,q)$ equals
$$
\sum_{k=0}^s  {\binom {d}{k}}_q q^{\binom {k}{2}} (1-p)^{q^{s-k}}
$$

(b)
The expected number $\mathrm{sp}(d,p,q)$ of subsets spanning $\mathrm{GF}(q)^d$
in a random subset $V(d,p,q)$ of vectors equals
$$
\sum_{k=0}^d {\binom {d}{k}}_q (-1)^k q^{\binom {k}{2}} (1 + p)^{q^{d-k}}.
$$

{\it Proof.}
Since the number of $s$-tuples distinguishing a matroid $M$ represented as a
multiset of vectors in $\mathrm{GF}(q)^d$ is the normalized
characteristic polynomial $\chi^{\dagger}(M;q^s)$ evaluated at $q^s,$
$$
\mathrm{D}(d,p,q,s) = \SC (\mathrm{GF}(q)^d;-p,\lambda).
$$ 
When $\lambda = q^s,$ the infinite product in the generating function telescopes into a 
finite product.  Thus,
$$
\sum_{d=0}^{\infty}  D(d,p,q,s) \frac {z^d}{[d!]_q}
=
\left[ \prod_{d=0}^{s-1} ( 1 + zq^d) \right]
\left[
\sum_{d=0}^{\infty}  (1-p)^{q^d} \frac {z^d}{[d!]_q}
\right].
$$
By an identity attributed to Euler or Cauchy (see, for example, \cite{GRota}, p.~254),
$$
\prod_{d=0}^{s-1} (1 - zq^d) =
\sum_{k=0}^{s} (-1)^k {\binom {d}{k}}_q q^{\binom {k}{2}} z^k.
$$
The formula now follows from replacing the product with the sum and expanding.

For part (b),  observe that $\SC (M|T; 1,0)$ is the number of subsets
$A$ in $T$ such that $r - \rk(A) = 0,$ that is, $A$ spans $M.$
Setting $\lambda = 0,$
$y_e = p_e,$ and $x_e = 1$ in Identity 7, we have
$$
\mathrm{sp}(d,p,q) = \SC (\mathrm{GF}(q)^d;p,0).
$$
From the generating function, in the form given in eqn (3), it follows that
$$
\sum_{d=0}^{\infty}  \mathrm{sp}(d,p,q) \frac {z^d}{[d!]_q}
=
\left[
\sum_{d=0}^{\infty}
(-1)^d q^{\binom {d}{2}} \frac {z^d}{[d!]_q}
\right]
\left[
\sum_{d=0}^{\infty}  (1+p)^{q^d} \frac {z^d}{[d!]_q}
\right].
$$
Expanding the product yields the formula.
This formula can also be obtained by M\"obius inversion on
the lattice $L(d,q).$
\qed

\vskip .20in \noindent
Identities 4 and 5 can also be given probabilistic interpretations.  If $M$ is a matroid on the set $E,$ we construct
a {\sl random contraction} $H$ by choosing a random subset $T$
by choosing each element $e$ independently and at random with probability $p_e$
and letting $H = M/T.$  A typical result is that the expected value of the
characteristic polynomial $\chi(H;\xi)$ of a random contraction is
$$
\underline{(1-p)}^E \mathbf{SC}
\left(M; \underline{\frac {p}{1-p}-1};\xi
\right).
$$

\vskip .33in \noindent
\begin{center}
4.  The motivating identities
\end{center}

We end by deriving the two identities which motivated this paper.
We begin by setting $x = -x,$ $y = -1,$ and $\lambda = 1,$ leaving $\xi$ unchanged in
Identity 3. Doing so, we obtain
the following convolutional identity, explicitly stated by Kook, Reiner, and Stanton
in
\cite{Kook} and implicit in the combinatorial construction of \'Etienne and Las Vergnas
in \cite{LaV}:
$$
\RR (M;x,\xi)
=
\sum_{T:\, T \subseteq E} \RR (M|T;x,-1)
\RR (M/T;-1,\xi).
\eqno(4)$$
We can restrict the range of summation to flats for the same reason as
in Identity 5.  Because $\RR (M/T;-1,\xi)$ is zero
if $M$ has isthmuses (or coloops), the range can be further restricted to
{\sl cyclic} flats, that is, flats with no isthmuses.

Next recall that
$$
\RR(M;x,-1) = \RR(M^{\perp};-1,x) = (-1)^{|E|-r}\chi(M^{\perp};-x),
$$
where $M^{\perp}$ is the orthogonal dual of $M.$  Thus, we obtain a
version of the ``KRSEV'' identity which expresses the rank generating polynomial
as a weighted-convolution of characteristic polynomials.

\vskip 0.2in\noindent
{\bf Identity 8.}
$$
(-1)^r \RR (M;-x,-\xi)
=
\sum_{X:\, X \in L^\circ(M)} (-1)^{|X|}  \chi ((M|X)^{\perp}; x)
\chi (M/X;\xi),
$$
where $L^\circ(M)$ is the lattice of cyclic flats of $M.$

\vskip 0.2in\noindent
Next, we set $x=1,$ $y=-1,$ leaving $\lambda$ and $\xi$ unchanged in Identity 1,
rederiving a multiplication identity in \cite{Kung}.

\vskip 0.2in\noindent
{\bf Identity 9.}
$$
\chi (M;\lambda\xi)
=
\sum_{X:\, X \in L(M)} \lambda^{r-\rk(X)}  \chi (M|X; \lambda)
\chi (M/X;\xi).
$$

Kook, Reiner, and Stanton \cite{Kook} have given an interpretation of Identity 8 for
graphs using pairs of colorings and flows.  One can easily adapt their interpretation
using the critical problem.  We will
just give the result, referring the reader to \cite{KungC}, Section 4.7, for
the necessary background.
Let $M$ be a rank-$r$ matroid on the set $E$ and suppose that $Q$ is an
$r \times |E|$ matrix over $\mathrm{GF}(s)$ representing $M$ and $R$ be an
$(|E|-r) \times |E|$ matrix over $\mathrm{GF}(t)$ representing the dual of $M.$
If $T \subseteq E,$ an {\sl $(s,t)$-duet $(u,v)$ with support $T$} is a pair of
row vectors, such that $u$ an $|E|$-dimensional vector
in the row space of $Q$ with all coordinates in
$T$ equal to $0$ and all coordinates in $E \backslash T$ not equal to $0,$ and
$v$ is a $|T|$-dimensional vector in the row space of the of  $R|T,$ the submatrix of $R$ consisting of all the columns labeled by $T,$ with all coordinates
not equal to $0.$   Then
$$
\RR (M;t,s)
=
(-1)^r  \sum_{(u,v)} (-1)^{|\mathrm{support}(T)|},
$$
where the sum ranges over all $(s,t)$-duets $(u,v).$

Interpretations of Identity 9 for graphs and representable matroids can be
found in \cite{Kung}.

\end{document}